\documentclass[12pt]{amsart}
\usepackage{fullpage,url,amssymb}

\DeclareFontEncoding{OT2}{}{} 

\usepackage{color}

\newcommand{\defi}[1]{\textsf{#1}} 

\newcommand{\PP}{{\mathbb P}}
\newcommand{\Q}{{\mathbb Q}}

\DeclareMathOperator{\Char}{char}

\newcommand{\tH}{{\operatorname{th}}}

\newcommand{\GL}{\operatorname{GL}}

\newcommand{\del}{\partial}

\newtheorem{theorem}{Theorem}[section]
\newtheorem{lemma}[theorem]{Lemma}

\theoremstyle{definition}

\theoremstyle{remark}
\newtheorem{remark}[theorem]{Remark}

\usepackage[
	backref,
	pdfauthor={Bjorn Poonen}, 
]{hyperref}
\usepackage[alphabetic,backrefs,lite]{amsrefs} 

\begin{document}

\title{Multivariable polynomial injections on rational numbers}
\author{Bjorn Poonen}
\thanks{This research was supported by NSF grant DMS-0841321.}
\address{Department of Mathematics, Massachusetts Institute of Technology, Cambridge, MA 02139-4307, USA}
\email{poonen@math.mit.edu}
\urladdr{http://math.mit.edu/~poonen}
\date{May 16, 2010}

\begin{abstract}
For each number field $k$,
the Bombieri-Lang conjecture for $k$-rational points
on surfaces of general type
implies the existence of a polynomial $f(x,y) \in k[x,y]$
inducing an injection $k \times k \to k$.
\end{abstract}

\maketitle

\section{Introduction}
\label{S:introduction}

Harvey Friedman asked whether there exists a polynomial 
$f(x,y) \in \Q[x,y]$ such that the induced map $\Q \times \Q \to \Q$
is injective.
Heuristics suggest that most sufficiently complicated polynomials
should do the trick.
Don Zagier has speculated that a polynomial as simple as $x^7+3y^7$ 
might already be an example.
But it seems very difficult to prove that {\em any} polynomial works.
Both Friedman's question and Zagier's speculation are at least a decade old 
(see~\cite{Cornelissen1999}*{Remarque~10}),
but it seems that there has been essentially no progress 
on the question so far.

Our theorem gives a positive answer conditional on 
a small part of a well-known conjecture.

\begin{theorem}
\label{T:main}
Let $k$ be a number field.
Suppose that there exists a homogeneous polynomial 
$F(x,y) \in k[x,y]$ 
such that the $k$-rational points on the surface $X$ in $\PP^3$
defined by $F(x,y)=F(z,w)$
are not Zariski dense in $X$.
Then there exists a polynomial $f(x,y) \in k[x,y]$
inducing an injection $k \times k \to k$.
\end{theorem}

\begin{remark}
If $F(x,y)$ is separable (or equivalently, squarefree) 
and homogeneous of degree at least $5$,
then $X$ is of general type.
So the hypothesis in Theorem~\ref{T:main}
would follow from the Bombieri-Lang conjecture
that $k$-rational points on a surface of general type are never Zariski dense.
\end{remark}

\begin{remark}
As the proof of Theorem~\ref{T:main} will show, 
if we have an algorithm for determining 
the Zariski closure of the set of $k$-rational points on each curve or surface
of general type, then we can construct $f(x,y)$ explicitly.
\end{remark}

\begin{remark}
To prove that a nonzero homogeneous polynomial $F(x,y)$ 
defines an injection $k \times k \to k$
is to prove that $X(k)$ is contained in the line $x-z=y-w=0$.
If $F$ is separable, then $X$ is a smooth projective hypersurface in $\PP^3$,
so it is simply connected.
But as far as we know, there is not a single simply connected 
smooth algebraic surface $X$ with $X(k) \not=\emptyset$
such that $X(k)$ is {\em known} to be not Zariski dense in $X$!
If one uses nonhomogeneous polynomials, one must instead
understand rational points on affine $3$-folds;
this seems unlikely to improve the situation.
All this suggests that Friedman's question cannot be answered unconditionally
without a major advance in arithmetic geometry.
\end{remark}

\begin{remark}
One cannot hope to answer the question using local methods alone.
More precisely, if $L$ is any local field of characteristic $0$,
and $f(x,y) \in L[x,y]$ is nonconstant,
then the induced map $L \times L \to L$ is not injective.
To prove this, choose a point $(x_0,y_0) \in L \times L$
where ${\del f}/{\del x}$ or ${\del f}/{\del y}$ is nonvanishing,
and let $c = f(x_0,y_0)$; then the affine curve $f(x,y)=c$ 
is smooth at $(x_0,y_0)$, so by the implicit function theorem
it contains infinitely many $L$-points, each of which has
the same image under $f$ as $(x_0,y_0)$.
\end{remark}

\begin{remark}
If $k$ is any imperfect field, then there exists a polynomial injection
$k \times k \to k$, by a construction that can be found in the proof
of Proposition~8 in \cite{Cornelissen1999}.
Namely, let $p=\Char k$, choose $t \in k - k^p$,
and use $f(x,y)=x^p+ty^p$.
This applies in particular to any global function field.
\end{remark}

\begin{remark}
The generalized $abc$-conjecture of \cite{Browkin-Brzezinski1994}
(more specifically, the $4$-variable analogue)
would imply that $f(x,y):=x^n+3y^n$ defines a polynomial injection
$\Q \times \Q \to \Q$ for sufficiently large odd integers $n$: 
this was observed in \cite{Cornelissen1999}*{Remarque~10}.
\end{remark}

\begin{remark}
For the function field $K$ of an irreducible curve 
over a base field $k$ of characteristic~$0$,
an analogue of the generalized $abc$-conjecture 
is known \cite{Mason1986}*{Lemma 2}.
This analogue can be used to show that for some $t \in K$ and $m \ge 1$,
the polynomial $f(x,y)=x^m+ty^m$ defines an injection,
under certain technical hypotheses.
These hypotheses can be satisfied when $k$ is a number field, for instance.
See \cite{Cornelissen1999}*{Proposition~8} for details
and for other related results.
\end{remark}

\section{Proof of theorem}

Let $k$, $F$, and $X$ be as in Theorem~\ref{T:main}.
Let $d=\deg F$.
Call a line in $\PP^3$ \defi{trivial} if it is given by
$x-\zeta z = y-\zeta w = 0$ for some $\zeta \in k$
with $\zeta^d=1$.
Each trivial line is contained in $X$.
Let $w$ be the number of roots of $1$ in $k$,
and let $p$ be a prime number such that $p>3$ and $p \nmid w$.
When we speak of the genus of a geometrically irreducible curve,
we mean the genus of its smooth projective model.
When we say that something holds for ``most'' elements of $k$
or of $k^n$, we mean that it holds outside a thin set
in the sense of~\cite{SerreMordellWeil1997}*{\S9.1}.
Such sets arise in the context of the Hilbert irreducibility theorem,
which shows that a finite union of thin sets cannot cover all of $k^n$.

\begin{lemma}
\label{L:inverse image}
Fix an integral closed subscheme $Z$ of $X$.
For most $\left(\begin{smallmatrix} a & b \\ c & d \end{smallmatrix}\right) \in \GL_2(k) \subset k^4$,
the inverse image $Y$ of $Z$ under the finite morphism
\begin{align*}
  \PP^3 &\to \PP^3 \\
  (x:y:z:w) &\mapsto (ax^p+by^p:cx^p+dy^p:az^p+bw^p:cz^p+dw^p)
\end{align*}
satisfies:
\begin{enumerate}
\item[(i)] If $\dim Z=0$, then $Y(k)=\emptyset$.
\item[(ii)] If $Z$ is a trivial line, then $Y(k)$ is contained in a trivial line.
\item[(iii)] If $Z$ is any other curve in $X$, then $Y(k)$ is finite.
\end{enumerate}
\end{lemma}

\begin{proof}
We can compute $Y$ in stages,
by first taking the {\em forward} image of $Z$ under the automorphism
\begin{align*}
  \PP^3 &\stackrel{\alpha}\to \PP^3 \\
  (x:y:z:w) &\mapsto (ax+by:cx+dy:az+bw:cz+dw)
\end{align*}
(technically, $\left(\begin{smallmatrix} a & b \\ c & d \end{smallmatrix}\right)$ here should be the inverse of what it was before, but this does not matter),
and then pulling back by 
\begin{align*}
  \PP^3 &\stackrel{\beta}\to \PP^3 \\
  (x:y:z:w) &\mapsto (x^p:y^p:z^p:w^p).
\end{align*}

(i) Here $\dim Z=0$.
If $Z$ is not a $k$-rational point, 
then $Z(k)=\emptyset$, so $Y(k)=\emptyset$.
If $Z$ is a $k$-rational point,
then for most
$\left(\begin{smallmatrix} a & b \\ c & d \end{smallmatrix}\right)$
the value of $(ax+by)/(cx+dy)$ on $Z$ is not a $p^\tH$ power in $k$,
so $Y(k)=\emptyset$.

(ii)
Here $Z$ is 
$x-\zeta z = y-\zeta w = 0$ for some $\zeta \in k$
with $\zeta^d=1$.
Then $\alpha(Z)=Z$, so $Y$ is $x^p-\zeta z^p = y^p -\zeta w^p = 0$.
By choice of $p$, 
the $p^\tH$-power map on $k$ is injective,
and moreover $\zeta=\eta^p$ for some $\eta \in k$ with $\eta^d=1$.
So all points in $Y(k)$ satisfy $x-\eta z = y-\eta w = 0$.

(iii)
Here $Z$ is an irreducible curve in $X$ that is not a trivial line.
If $Z$ is geometrically {\em reducible},
then $Z(k)$ is not Zariski dense in $Z$, so $Z(k)$ is finite,
and $Y(k)$ is finite too.
So assume that $Z$ is geometrically irreducible.

If $y=0$ on $Z$ or if $x/y$ defines a {\em constant} rational function on $Z$,
then as in (i),
for most $\left(\begin{smallmatrix} a & b \\ c & d \end{smallmatrix}\right)$
the value of $x/y$ on $\alpha(Z)$ on $Z$ is not a $p^\tH$ power in $k$,
so $Y$ has no $k$-rational points except possibly those where $x=y=0$,
so $Y(k)$ is finite.

Suppose that $x/y$ 
defines a rational function of degree $m>1$ on $Z$.
By Bertini's theorem (\cite{Hartshorne1977}*{Corollary~III.10.9}),
$ax+by$ has distinct zeros on the normalization $Z'$ of $Z$,
outside the base locus of the linear system given by $\langle x,y \rangle$,
for most $a$ and $b$.
The same applies to $cx+dy$, so for most 
$\left(\begin{smallmatrix} a & b \\ c & d \end{smallmatrix}\right)$,
the rational function $(ax+by)/(cx+dy)$ on $Z'$ has $m$ simple zeros
and $m$ simple poles on $Z'$.
Adjoining the $p^\tH$ root of this function to the function field of $Z'$
yields the function field of a geometrically irreducible curve $C$ 
of genus greater than $1$, by the Hurwitz formula.
By~\cite{Faltings1983}, $C(k)$ is finite.
Since $Y$ admits a dominant rational map to $C$, the set $Y(k)$ is finite too.

Thus we may assume that $x/y$ is of degree $1$ on $Z$;
in particular, $Z$ is a rational curve.
Similarly, we may assume that $z/w$ is of degree $1$ on $Z$.
If the rational functions $x/y$ and $z/w$ on $Z$ were different,
then for most $\left(\begin{smallmatrix} a & b \\ c & d \end{smallmatrix}\right)$,
the supports of the divisors of $(ax+by)/(cx+dy)$ and $(az+bw)/(cz+dw)$
on the normalization of $Z$ would not coincide.
Adjoining the $p^\tH$ roots of these functions would lead to a geometrically
irreducible curve of genus greater than $1$, by the Hurwitz formula again.
So $Y(k)$ would be finite as before.

Thus we may assume that $x/y=z/w$ as rational functions on $Z$.
So on $Z$, we have
\[
	x^d F(x,y) = x^d F(z,w) = F(xz,xw) = F(xz,yz) = z^d F(x,y).
\]
But $F(x,y)$ does not vanish on $Z$ (since $x/y$ is nonconstant),
so $x^d-z^d$ vanishes on $Z$.
Since $Z$ is geometrically irreducible, $x-\zeta z$ vanishes on $Z$
for some $\zeta \in k$ with $\zeta^d=1$.
But $x/y=z/w$ on $Z$, so $y-\zeta w$ vanishes on $Z$ too.
Thus $Z$ is a trivial line, a contradiction.
\end{proof}

Let $W$ be the Zariski closure of $X(k)$.
By assumption, $\dim W \le 1$.
Applying Lemma~\ref{L:inverse image} to each irreducible component of $W$
shows that by replacing $F(x,y)$ 
with $F(ax^p+by^p,cx^p+dy^p)$ for suitable 
$\left(\begin{smallmatrix} a & b \\ c & d \end{smallmatrix}\right)$,
we may reduce to the case that $W(k)$ contains at most finitely many
points outside the trivial lines.
Repeating this construction lets us reduce to the case 
that $W(k)$ is contained in the trivial lines.

\begin{lemma}
\label{L:finitely many exceptions}
For $G(x,y):=F(x^p+1,y^p+1)$ for the new $F$ as above,
the equation $G(x,y)=G(z,w)$ has only finitely many solutions over $k$
with $(x,y) \ne (z,w)$.
\end{lemma}

\begin{proof}
Suppose that $x,y,z,w \in k$ are such that $G(x,y)=G(z,w)$.
Then there exists $\zeta \in k$ with $\zeta^d=1$
such that $x^p+1 = \zeta(z^p+1)$ and $y^p+1 = \zeta(w^p+1)$.
If $\zeta=1$, then this implies $(x,y)=(z,w)$, by the choice of $p$.
If $\zeta \ne 1$,
then $x^p+1 = \zeta(z^p+1)$ defines a geometrically irreducible curve 
whose projective closure is smooth, and hence of genus 
$p(p-1)/2 > 1$, so by~\cite{Faltings1983} 
it has at most finitely many solutions in $k$.
The same applies to $y^p+1=\zeta(w^p+1)$,
and there are only finitely many $\zeta$.
\end{proof}

\begin{lemma}
\label{L:eliminating exceptions}
If $G(x,y) \in k[x,y]$ is such that 
the equation $G(x,y)=G(z,w)$ has only finitely many solutions over $k$
with $(x,y) \ne (z,w)$,
then for most $(a,b) \in k^2$
the polynomial $f(x,y):=G(a x^p + b, a y^p + b)$
defines an injection $k \times k \to k$.
\end{lemma}

\begin{proof}
By choice of $p$, if $a \ne 0$ and $b$ is arbitrary,
the polynomial $ax^p+b$ defines an injection $k \to k$.
Also, for most $(a,b)$, 
the coordinates of the finitely many exceptional solutions to $G(x,y)=G(z,w)$ 
are not in the range of this injection.
\end{proof}

Theorem~\ref{T:main} follows from 
Lemmas \ref{L:finitely many exceptions} and~\ref{L:eliminating exceptions}.

\section*{Acknowledgements} 

The idea for this article arose during a discussion session at
the Hausdorff Institute in Bonn.
I thank Burt Totaro for a comment on a first draft of this article,
and I thank the referee for a few suggestions.

\end{document}